\documentclass[12pt]{iopart}

\usepackage{amssymb,amsmath,amsfonts,amssymb}
\usepackage{graphics,graphicx,color}

\def\@abssec#1{\vspace{.05in}\footnotesize \parindent .2in
{\bf #1. }\ignorespaces}

\graphicspath{{/EPSF/}{Figures/}}
\DeclareGraphicsExtensions{.eps}

\newtheorem{theorem}{Theorem}[section]
\newtheorem{lemma}[theorem]{Lemma}
\newtheorem{proposition}[theorem]{Proposition}
\newtheorem{corollary}[theorem]{Corollary}

\newtheorem{hypothesis}[theorem]{Hypothesis}
\def \Rm {\mathbb R}
\def \Nm {\mathbb N}
\def \Cm {\mathbb C}

\def\S{\mathbb S}

\def\pa{\partial}

\def\d{{\rm dist}}
\newcommand{\eps}{\varepsilon}

\newcommand{\dsum}{\displaystyle\sum}
\newcommand{\dint}{\displaystyle\int}
\newcommand{\pdr}[2]{\dfrac{\partial{#1}}{\partial{#2}}}
\newcommand{\pdrr}[2]{\dfrac{\partial^2{#1}}{\partial{#2}^2}}

\newcommand{\aver}[1]{\langle {#1} \rangle}

\renewcommand{\d}{\mathfrak d}
\renewcommand{\k}{\mathfrak k}  
\newcommand{\vp}{\varphi}
\newcommand{\dvp}{\dot\vp}

\newcommand{\cout}[1]{}

 \renewcommand{\arraystretch}{1.5}

\begin{document}
\title{Inverse Diffusion Theory of Photoacoustics}

\author{Guillaume Bal}
\address{Department of Applied Physics \& Applied Mathematics,
   Columbia University, New York, NY 10027}
 \ead{gb2030@columbia.edu}

\author{Gunther Uhlmann}
\address{Department of Mathematics,
  University of Washington, Seattle, WA 98195}
 \ead{gunther@math.washington.edu}



\begin{abstract}
  This paper analyzes the reconstruction of diffusion and absorption
  parameters in an elliptic equation from knowledge of internal data.
  In the application of photoacoustics, the internal data are the
  amount of thermal energy deposited by high frequency radiation
  propagating inside a domain of interest. These data are obtained by
  solving an inverse wave equation, which is well-studied in the
  literature.
  
  We show that knowledge of two internal data based on well-chosen
  boundary conditions uniquely determines two constitutive parameters
  in diffusion and Schr\"odinger equations. Stability of the
  reconstruction is guaranteed under additional geometric constraints
  of strict convexity. No geometric constraints are necessary when
  $2n$ internal data for well-chosen boundary conditions are
  available, where $n$ is spatial dimension. The set of well-chosen
  boundary conditions is characterized in terms of appropriate complex
  geometrical optics (CGO) solutions.
  
\end{abstract}
 

\renewcommand{\thefootnote}{\fnsymbol{footnote}}
\renewcommand{\thefootnote}{\arabic{footnote}}

\renewcommand{\arraystretch}{1.1}




\section{Introduction}
\label{sec:intro}

Photoacoustic tomography (PAT) is a recent hybrid medical imaging
modality that combines the high resolution of acoustic waves with the
large contrast of optical waves. When a body is exposed to short
pulse radiation, typically emitted in the near infra-red region in
PAT, it absorbs energy and expands thermo-elastically by a very small
amount; this is the photoacoustic effect. Such an expansion is
sufficient to emit acoustic pulses, which travel back to the boundary
of the domain of interest where they are measured.

A first step in PAT is therefore to reconstruct the amount of
deposited energy from time-dependent boundary measurement of acoustic
signals. Acoustic signals propagate in fairly homogeneous domains as
the sound speed varies very little from one tissue to the next.  The
reconstruction of the amount of deposited energy is therefore quite
accurate in many practical settings. For references on the practical
and theoretical aspects of PAT, we refer the reader to e.g.
\cite{ABJK-SR-09,CAB-IP-07,CAB-JOSA-09,CLB-SPIE-09,FPR-JMA-04,FR-CRC-09,FSS-PRE-07,HSBP-IP-04,HSS-MMAS-05,HKN-IP-08,KK-EJAM-08,PS-IP-07,SU-IP-09,XW-RSI-06,XWKA-CRC-09}.

Once the amount of deposited energy has been reconstructed, a second
step consists of inferring the optical properties in the body. Although
this second step is less studied mathematically, it has received
significant attention in the biomedical literature; see e.g.
\cite{CAB-IP-07,CAB-JOSA-09,CLB-SPIE-09}. The reconstruction of the
optical parameters is very useful in practice because the attenuation
properties of healthy and unhealthy tissues are extremely different
\cite{CLB-SPIE-09,XW-RSI-06}. The combination of high spatial
resolution of the acoustic signals and large contrast of the optical
parameters allows us e.g. to image blood vessels at high resolution,
which is in turn important for assessing the status of cancerous
tissues.

Near infra-red light is best modeled by radiative transfer equations.
What we can reconstruct on the optical coefficients, seen as
constitutive parameters in the radiative transfer equation, from
measurements of heat deposition for all possible illuminations of the
domain has recently been analyzed in \cite{BJJ-IP-09}. The latter
paper considers the continuous-illumination setting. In other words,
heat deposition is measured for all possible incoming radiation
condition, which can be controlled at the boundary of the domain.

In this paper, we consider the diffusion approximation to radiative
transfer, which is accurate for radiation propagation in highly
scattering media. Such an approximation is typically valid for
propagations of radiation over one centimeter or more
\cite{arridge99,B-IP-09}.  In such a simplified setting, the unknown
optical parameters are the spatially varying diffusion and attenuation
coefficients.  Provided that the diffusion coefficient is known at the
domain's boundary, we show that these two coefficients are uniquely
determined by two well-chosen illuminations at the boundary of the
domain. An explicit reconstruction procedure is presented, which
involves solving a first-order equation with a vector field explicitly
constructed from the internal data (the amount of deposited energy).

The stability of the reconstruction of the optical coefficients from
two internal data is established under geometric conditions of strict
convexity on the domain of interest. In the presence of $2n$
well-chosen boundary conditions, where $n$ is spatial dimension, such
geometric constraints can be removed and constructions are possible on
essentially arbitrary domains of interest with well-defined boundary.

Mathematically, the inverse problem is an inverse diffusion problem
with two internal measurements. By using the standard Liouville change
of variables, the diffusion equation is replaced by a Schr\"odinger
equation with unknown potential and with internal measurements
involving a second unknown function. By adapting the theory of complex
geometrical optics solutions to this setting, we are able to obtain
uniqueness and stability results for the inverse Schr\"odinger
problem. By means of the inverse Liouville change of variables, we are
able to conclude on the uniqueness and stability of the reconstruction
of the optical parameters for the photoacoustic tomography problem in
the diffusive regime.

The rest of the paper is structured as follows. Section \ref{sec:main}
presents the photoacoustic tomography problem and the main results we
obtain in this paper. The inverse Schr\"odinger problem and the
explicit reconstruction algorithms are addressed in section
\ref{sec:Sch}. How the latter results are used to solve the 
inverse diffusion problem is explained in section \ref{sec:idid}.

\section{Photoacoustic tomography and main results}
\label{sec:main}

The propagation of radiation in scattering media is modeled by the 
following diffusion equation
\begin{equation}
  \label{eq:timediffusion}
  \begin{array}{l}
    \dfrac{1}{c} \pdr{}t u
   -\nabla\cdot D(x) \nabla u + \sigma_a(x) u =0 ,\qquad x\in X\subset\Rm^n,
    \,\, t\geq0\\
    u = g \qquad x\ \in \partial X,\,\, t\geq0
  \end{array}
\end{equation}
where $X$ is an open, bounded, connected domain in $\Rm^n$ with $C^1$
boundary $\pa X$ (embedded in $\Rm^n$), where $n$ spatial dimension;
$c$ is light speed in tissues; $D(x)$ is a (scalar) diffusion
coefficient; and $\sigma_a(x)$ is an attenuation coefficient. The
source of incoming radiation is prescribed by $g(t,x)$ on the boundary
$\pa X$ and is assumed to be a very short pulse supported on an
interval of time $(0,\eta)$ with $c\eta$ of order $O(1)$. The amount
of energy deposited is proportional to attenuation and is given by
\begin{displaymath}
  H(t,x) = \sigma_a(x) u(t,x).
\end{displaymath}
A thermal expansion (assumed to be proportional to $H$) results and
emits acoustic waves. Such waves are modeled by
\begin{equation}
  \label{eq:wave}
  \dfrac{1}{c_s^2(x)} \pdrr pt - \Delta p = \beta \pdr{}t H(t,x),
\end{equation}
with $c_s$ the sound speed and $\beta$ a coupling coefficient assumed
to be constant and known. The pressure (potential) $p(t,x)$ is then
measured on $\pa X$ as a function of time. Using as in
\cite{BJJ-IP-09} the difference of time scales $c_s\ll c$, which shows
that radiation propagation occurs at a much faster time scale than
acoustic wave propagation, we can show that
\begin{displaymath}
  H(t,x) \sim H_0(x) \delta_0(t),\qquad H_0(x) = \sigma_a(x)\dint_{\Rm^+}
    u(t,x)dt.
\end{displaymath}
We now have a well-posed inverse wave problem, where measurement of
$n-$dimensional information $p(t,x)$ for $t>0$ and $x\in\pa X$ allows
us to reconstruct the $n-$dimensional spatial map $H_0(x)$.  We assume
this step, which has been analyzed in great detail in the literature,
completed; see e.g.
\cite{FPR-JMA-04,HKN-IP-08,KK-EJAM-08,SU-IP-09,XWKA-CRC-09}.  

Denoting by $g(x)=\int_{\Rm^+}g(t,x)dt$ and
$u(x)=\int_{\Rm^+}u(t,x)dt$, we thus observe that the photoacoustic
problem in the diffusive regime amounts to reconstructing
$(D(x),\sigma_a(x))$ from knowledge of 
\begin{equation}
  \label{eq:internmeas}
  H_j(x) = \sigma_a(x) u_j(x),\qquad 1\leq j\leq J,
\end{equation}
for $J\in\Nm^*$ illumination maps $g_j(x)$, where $u_j$ is the 
solution to the steady-state equation
\begin{equation}
  \label{eq:diffusion}
  \begin{array}{l}
   -\nabla\cdot D(x) \nabla u_j + \sigma_a(x) u_j =0 ,
   \qquad x\in X\subset\Rm^n,\\
    u_j = g_j \qquad x \in \partial X.
  \end{array}
\end{equation}
The set of internal data is given by:
\begin{equation}
  \label{eq:dj}
  d=(d_j)_{1\leq j\leq J}, \qquad d_j(x) = \sigma_a(x) u_j(x).
\end{equation}

The main problem of interest in this paper is the uniqueness and the
stability of the reconstruction for the Inverse Diffusion problem
with Internal Data:
\begin{center}
  {\bf ISID:} Reconstruction of $(D(x),\sigma_a(x))$ from knowledge of
  $d=(d_j)_{1\leq j\leq J}$ on $X$ for a {\em fixed} collection of
  illuminations $g=(g_j)_{1\leq j\leq J}$ prescribed on $\pa X$.
\end{center}
We start with the case of two real-valued measurements $J=2$. The
reconstructions are based on the construction of vector fields that
are well defined only when the optical coefficients are sufficiently
smooth. More precisely, let $k\geq1$ and assume that
\begin{displaymath}
  \sqrt D\in Y=H^{\frac n2+k+2+\eps}(X) \subset C^{k+2}(\bar X),\quad
  \sigma_a \in C^{k+1}(\bar X),\qquad \eps>0.
\end{displaymath}
Let then $\k\in\Rm^n$ be a constant vector with $|\k|$ sufficiently
large.

Let $g=(g_1,g_2)$ be a given illumination and $\d=d_1+id_2=
\sigma_a(x)(u_1+iu_2)$ the corresponding internal data, where the
real-valued solutions $u_{j}$ solve \eqref{eq:diffusion} with boundary
conditions $g_{j}$ for $j=1,2$. We introduce the vector field
and scalar quantity
\begin{equation}
  \label{eq:betagammak}
  \begin{array}{rclrcl}
    \beta &=& \dfrac{e^{-2\k\cdot x}}{2|\k|} \Im\big(\d\nabla \bar{\d}
    -\bar {\d}\nabla \d\big),\qquad &
    \gamma &=&  \dfrac{e^{-2\k\cdot x}}{4|\k|} \Im
       \big(\bar{\d}\Delta \d-\d\Delta{\bar \d}\big).
  \end{array}
\end{equation}
Here, $\bar \d=d_1-id_2$ is the complex conjugate of $\d$.  Then we will
see in section \ref{sec:Sch} that for $\mu=D^{-\frac12}\sigma_a$, we
have
\begin{equation}
  \label{eq:mu000}
  \beta\cdot\nabla \mu + \gamma \mu =0.
\end{equation}
Since $\mu$ is known on $\pa X$ as $D$ and $\sigma_a=d_j/g_j$ are
known on $\pa X$, the above equation is a well-posed equation for
$\mu$.  When the integral curves of the vector field $\beta$ map any
point inside $x$ to a point $x_0(x)\in \pa X$, then \eqref{eq:mu000}
uniquely characterizes $\mu$ inside the domain.

Let $v=\sqrt D \Re u$ and define $q=-v^{-1}\Delta v$, which we will
show is well-defined. Then it turns out that 
\begin{displaymath}
  -\Delta \sqrt D - q \sqrt D = \mu,\qquad x\in X.
\end{displaymath}
This provides an elliptic equation for $\sqrt D$ that admits a unique
solution since $\sqrt D$ is known on $\pa X$.

The principal difficulty is to ensure that $\beta$ is constructed in
such a way that \eqref{eq:mu000} allows for a unique solution $\mu$.
We will show that there exists an open set of boundary conditions $g$
on $\pa X$ that allows us to do so. We define the set of coefficients
$(D,\sigma_a)\in \mathcal M$ as:
\begin{equation}
  \label{eq:M000}
  \mathcal M = \{(D,\sigma_a)\mbox{ such that }
  (\sqrt D,\sigma_a)\in Y\times C^{k+1}(\bar X), 
    \,\,
   \|\sqrt D\|_{Y} + \|\sigma_a\|_{C^{k+1}(\bar X)} \leq M \}.
\end{equation}
  
The main results for the Inverse Diffusion problem with Internal Data
(IDID) are then as follows. We start with a uniqueness results:
\begin{theorem}
  \label{thm:uniq000}
  Let $X$ be an open, bounded, domain with $C^2$ boundary $\pa X$.
  Assume that $(D,\sigma_a)$ and $({\tilde D},\tilde\sigma_a)$ are in
  $\mathcal M$ with $D_{|\pa X}=\tilde D_{|\pa X}$. Let $d$ and
  $\tilde d$ be the internal data in \eqref{eq:dj} for the
  coefficients $(D,\sigma_a)$ and $(\tilde D,\tilde\sigma_a)$,
  respectively and with boundary conditions $(g_j)_{j=1,2}$.  Then
  there is an open set of illuminations $g\in (C^{1,\alpha}(\pa X))^2$
  for some $\alpha>\frac12$ such that if $d=\tilde d$, then
  $(D,\sigma_a)=(\tilde D,\tilde\sigma_a)$.
\end{theorem}

The above result shows uniqueness of the reconstruction of the optical
coefficients but does not imply stability. Without additional
geometric information about $\pa X$, the above procedure may yield
unstable reconstructions. However, provided that $2n$ well-chosen
measurements $d_j$ for $1\leq j\leq J=2n$ are available, then
$n$ vector fields similar to $\beta$ above may be constructed in order
to form locally a basis of vectors in $\Rm^n$. In such a setting, the 
following stability result holds. 
\begin{theorem}
  \label{thm:stab1000}
  Let $k\geq2$ and let $X$ be an arbitrary bounded domain with
  boundary $\pa X$ of class $C^{k+1}$.  Assume that $(D,\sigma_a)$ and
  $({\tilde D},\tilde\sigma_a)$ are in $\mathcal M$ with
  $D_{|\pa X}=\tilde D_{|\pa X}$. Let $d=(d_1,\ldots, d_{2n})$ and
  $\tilde d=(\tilde d_1,\ldots,\tilde d_{2n})$ be the internal data
  constructed in \eqref{eq:dj} for the coefficients $(D,\sigma_a)$ and
  $(\tilde D,\tilde\sigma_a)$, respectively and with boundary
  conditions $g=(g_j)_{1\leq j\leq 2n}$.  Then there is an open set of
  illuminations $g\in (C^{k,\alpha}(\pa X))^{2n}$ and a constant $C$
  such that
  \begin{equation}
    \label{eq:stab1000}
     \|D-\tilde D\|_{C^{k}(X)} + \|\sigma_a-\tilde\sigma_a\|_{C^{k}(X)}
   \leq C \|d-\tilde d\|_{(C^{k+1}(X))^{2n}}.
  \end{equation}
\end{theorem}

Stability is therefore ensured without geometric constraints provided
that enough measurements are available. When only two measurements are
available and are of the form described in Theorem \ref{thm:uniq000},
we can still get a stability result under the following geometric
hypothesis:
\begin{hypothesis}
  \label{hyp:R0} There exists $R<\infty$  such that for each $x_0\in \pa X$, 
  which we assume is of class $C^2$, we have $X\subset B_{x_0}(R)$
  where $B_{x_0}(R)$ is a ball of radius $R$ that is tangent to $\pa
  X$ at $x_0\in\pa X$.
\end{hypothesis}
Then we can show the following result.
\begin{theorem}
  \label{thm:stab2000}
  Let $k\geq3$.  Let $X$ satisfy Hypothesis \ref{hyp:R0} with boundary
  $\pa X$ of class $C^{k+1}$.  Assume that $(D,\sigma_a)$ and
  $({\tilde D},\tilde\sigma_a)$ are in $\mathcal M$ with $D_{|\pa
    X}=\tilde D_{|\pa X}$. Let $d$ and $\tilde d$ be internal data as
  above for the coefficients $(D,\sigma_a)$ and $(\tilde
  D,\tilde\sigma_a)$, respectively and with boundary conditions
  $g=(g_j)_{j=1,2}$.  Then there is an open set of illuminations $g\in
  (C^{k,\alpha}(\pa X))^2$ and a constant $C$ such that
  \begin{equation}
    \label{eq:stab2000}
     \|D-\tilde D\|_{C^{k-1}(X)} + \|\sigma_a-\tilde\sigma_a\|_{C^{k-1}(X)}
   \leq C \|d-\tilde d\|_{(C^{k}(X))^2}.
  \end{equation}
\end{theorem}

The above three theorems are proved in the following two sections.
They show that the internal data $d$ for illuminations $g$ with
$J\geq2$ provide stable reconstructions of the optical coefficients
under the assumption that the illuminations are well-chosen. How these
illuminations are chosen will become more explicit in the next two
sections. The characterization of the open set of illuminations is
however not very precise. The main features of the result are as
follows. For coefficients in $\mathcal M$, there is a minimal value of
$|\k|$ that ensures that we can construct a vector field $\beta$ with
Property {\bf P}, which means its integral curves map any point in $X$
to an point in $\pa X$. Such a vector field is constructed by means of
complex geometric optics solutions, which depend on the unknown
optical parameters. The illuminations must then be chosen sufficiently
close to the trace on $\pa X$ of the above vector field to ensure that
they generate another vector field with Property {\bf P}.  Closedness
is therefore not characterized in very explicit means. It remains an
interesting question to obtain a priori constraints on the illumination
that will ensure that the resulting vector field satisfies Property
{\bf P}.

In the next section, we consider a similar problem for the
Schr\"odinger equation. How the latter results are used to prove the
theorems stated above is described in section \ref{sec:idid}.

\section{Inverse Schr\"odinger with Internal Data}
\label{sec:Sch}

Let $X$ be an open, bounded, connected, domain in $\Rm^n$, where $n$
is spatial dimension, with smooth boundary $\pa X$. We consider the
Schr\"odinger equations
\begin{equation}
  \label{eq:Sch}
  \begin{array}{ll}
    \Delta u_j + q u_j =0\qquad & X \\
    u_j=g_j & \partial X,
  \end{array}
\end{equation}
for $1\leq j\leq J$. Here, $J\in\Nm^*$ is the number of illuminations
and $q$ is an unknown potential. We assume that the homogeneous
problem with $g_j=0$ admits the unique solution $u\equiv0$ so that
$\lambda=0$ is not in the spectrum of $\Delta+q$.  We assume that $q$
on $X$ is the restriction to $X$ of a function $\tilde q$ compactly
supported on $\Rm^n$ and such that $\tilde q\in H^{\frac
  n2+k+\eps}(\Rm^n)$ with $\eps>0$ for $k\geq1$. Moreover we assume
that the extension is chosen so that
\begin{equation}
  \label{eq:bdext}
  \|q\|_{H^{\frac n2+k+\eps}(X)} \leq C(X,k,n)
  \|\tilde q\|_{H^{\frac n2+k+\eps}(\Rm^n)},
\end{equation}
for some constant $C(X,k,n)$ independent of $q$. That such a constant
exists may be found e.g. in \cite[Chapter VI, Theorem 5]{S-PUP-70}.

We assume that $g_j\in C^{k,\alpha}(X)$ with $\alpha>\frac12$ and $\pa
X$ is of class $C^{k+1}$ so that \eqref{eq:Sch} admits a unique
solution $u_j\in C^{k+1}(X)$ \cite[Theorem 6.19]{gt1}.  The internal
data are of the form
\begin{equation}
  \label{eq:meas}
  d_j(x) = \mu(x) u_j(x), \qquad X, \quad 1\leq j\leq J.
\end{equation}
Here $\mu\in C^{k+1}(\bar X)$ verifies $0<\mu_0\leq \mu(x)\leq \mu_0^{-1}$ for
a.e. $x\in X$. 

The {\em inverse Schr\"odinger problem with internal data} (ISID)
consists of reconstructing $(q,\mu)$ in $X$ from knowledge of
$d=(d_1,\ldots,d_J)\in \big(C^{k+1}(X)\big)^J$ for a given
illumination $g=(g_j)_{1\leq j\leq J}$. We will mostly be concerned
with the case $J=2$ and $J=2n$ with $g_j$, and hence $d_j$ real-valued
measurements.

%
%
%
%
%
\subsection{Complex Geometrical Optics Solutions}
\label{sec:CGO}

The analysis of ISID carried out in this paper is based on the
construction of complex geometrical optics (CGO) solutions. When
$q=0$, CGOs are harmonic solutions of the form $e^{\rho\cdot x}$ for
$\rho\in\Cm^n$ such that $\rho\cdot\rho=0$. When $q\not\equiv0$, CGOs
are solutions of the following problem
\begin{equation}
  \label{eq:CGO0}
  \Delta u_\rho + qu_\rho =0 ,\qquad u_\rho \sim e^{\rho\cdot x} \mbox{ as } 
    |x|\to\infty.
\end{equation}
More precisely, we say that $u_\rho$ is a solution of the above equation
with $\rho\cdot\rho=0$ and the proper behavior at infinity when it is
written as
\begin{equation}
  \label{eq:decurho}
  u_\rho(x) = e^{\rho\cdot x} \big(1+\psi_\rho(x)\big),
\end{equation}
for $\psi_\rho\in L^2_\delta$ a weak solution of
\begin{equation}
  \label{eq:psi}
  \Delta \psi_\rho + 2\rho\cdot\nabla \psi_\rho = -q(1+\psi_\rho).
\end{equation}
The space $L^2_\delta$ for $\delta\in\Rm$ is defined as the completion 
of $C^\infty_0(\Rm^n)$ with respect to the norm $\|\cdot\|_{L^2_\delta}$
defined as
\begin{equation}
  \label{eq:normdelta}
  \|u\|_{L^2_\delta}=\Big(\dint_{\Rm^n}
     \aver{x}^{2\delta}|u|^2 dx\Big)^{\frac12}, \qquad
   \aver{x}=(1+|x|^2)^{\frac12}.
\end{equation}
Let $-1<\delta<0$ and $q\in L^2_{\delta+1}$ and $\aver{x}q\in L^\infty$.
One of the main results in \cite{Syl-Uhl-87} is that there exists 
$\eta=\eta(\delta)$ such that the above problem admits a 
unique solution with $\psi\in L^2_\delta$ provided that 
\begin{displaymath}
  \|\aver{x}q\|_{L^\infty} +1 \leq \eta |\rho|.
\end{displaymath}
Moreover, $\|\psi\|_{L^2_\delta}\leq C|\rho|^{-1}\|q\|_{L^2_{\delta+1}}$
for some $C=C(\delta)$. In the analysis of ISID, we need smoother 
CGOs than what was recalled above. We introduce the spaces 
$H^s_\delta$ for $s\geq0$ as the  completion 
of $C^\infty_0(\Rm^n)$ with respect to the norm $\|\cdot\|_{H^s_\delta}$
defined as
\begin{equation}
  \label{eq:normHsdelta}
  \|u\|_{H^s_\delta}=\Big(\dint_{\Rm^n}
     \aver{x}^{2\delta}|(I-\Delta)^{\frac s2}u|^2 dx\Big)^{\frac12}.
\end{equation}
Here $(I-\Delta)^{\frac s2}u$ is defined as the inverse Fourier
transform of $\aver{\xi}^s\hat u(\xi)$, where $\hat u(\xi)$ is the Fourier
transform of $u(x)$. Then we have the following 
\begin{proposition}
  \label{prop:regul}
  Let $-1<\delta<0$ and $k\in\Nm^*$. Let $q\in H^{\frac n2+k+\eps}_{1}$ and
  hence in $ H^{\frac n2+k+\eps}_{\delta+1}$ and $\rho$ be such that
  \begin{equation} \label{eq:constq}
    \|q\|_{H_1^{\frac n2+k+\eps}} +1 \leq \eta |\rho|.
  \end{equation}
  Then $\psi_\rho$ the unique solution to \eqref{eq:psi} belongs
  to $H^{\frac n2+k+\eps}_\delta$ and 
  \begin{equation}
    \label{eq:cts}
    |\rho| \|\psi_\rho\|_{H^{\frac n2+k+\eps}_\delta}
   \leq C  \|q\|_{H^{\frac n2+k+\eps}_{\delta+1}}, 
  \end{equation}
  for a constant $C$ that depends on $\delta$ and $\eta$.
\end{proposition}
\begin{proof}
  We recall \cite{Syl-Uhl-87} that for $|\rho|\geq c>0$ and $f \in
  L^2_{\delta+1}$ with $-1<\delta<0$, the equation
  \begin{equation}
    \label{eq:inhom}
    (\Delta+2\rho\cdot\nabla) \psi =f 
  \end{equation}
  admits a unique weak solution $\psi\in L^2_{\delta}$ with
  \begin{displaymath}
    \|\psi\|_{L^2_\delta} \leq C(\delta,c) |\rho|^{-1} \|f\|_{L^2_{\delta+1}}.
  \end{displaymath}
  Now since $(\Delta+2\rho\cdot\nabla)$ and $(I-\Delta)^s$ are
  constant coefficient operators and hence commute, we deduce that 
  when $f\in H^s_{\delta+1}$ for any $s>0$, then
  \begin{equation} \label{eq:regctcoef}
     \|\psi\|_{H^s_\delta} \leq C(\delta,c) |\rho|^{-1} \|f\|_{H^s_{\delta+1}}.
  \end{equation}
  The solution to \eqref{eq:psi} is known to admit the decomposition
  \begin{displaymath}
    \psi_\rho = \dsum_{j=0}^\infty \psi_j, \qquad
       \Delta \psi_j + 2\rho\cdot\nabla \psi_j = -q\psi_{j-1},
  \end{displaymath}
  with $\psi_{-1}=1$. Let $s=\frac n2+k+\eps$. Assume $q\psi_{j-1}\in
  H^s_{\delta+1}$ which is true for $j=0$ by assumption on $q$.  Then
  $\psi_j\in H^s_{\delta}$. Since $H^{s}$ is an algebra, we want to
  prove that
  \begin{equation} \label{eq:ineqCGO}
    \|q\psi_j\|_{H^s_{\delta+1}} \leq \|q\|_{H^s_1}
   \|\psi_j\|_{H^s_\delta}.
  \end{equation}
  Indeed, decompose $\Rm^n$ into cubes. On each cube $B$,
  $\aver{x}^{2s}$ is more or less constant up to a $C^{\pm 2s}$. Now
  $H^s(B)$ is an algebra so that $\|qu\|_{H^s(B)}\leq
  \|q\|_{H^s(B)}\|u\|_{H^s(B)}$.  Since $\aver{x}$ is more or less
  constant and equal to $\aver{x_B}$,
  \begin{displaymath}
      \aver{x_B}^{2\delta+2}\|qu\|_{H^s(B)}^2 \leq C
      \|\aver{x}(I-\Delta)^{\frac s2}q\|^2_{L^2(B)}
      \|\aver{x}^{\delta}(I-\Delta)^{\frac s2}u\|^2_{L^2(B)}
  \end{displaymath}
  It remains to sum over all the cubes $B$ to get the result. When the
  size of the cubes tends to $0$, the constant $C$ tends to $1$, which
  yields \eqref{eq:ineqCGO}.
  This and \eqref{eq:regctcoef} show that 
  \begin{displaymath}
    \|\psi_j\|_{H^s_\delta} \leq C |\rho|^{-1} \|q\|_{H^{\frac s2+1+\eps}_1}
    \|\psi_{j-1}\|_{H^s_\delta}.
  \end{displaymath}
  By selecting $\eta$ such that $C |\rho|^{-1} \|q\|_{H^{\frac
      s2+k+\eps}_1}<\frac12$, we obtain 
  \begin{displaymath}
    \|\psi_j\|_{H^s_\delta} \leq \dfrac{1}{2^j} C |\rho|^{-1}
    \|q\|_{H^s_{\delta+1}}.
  \end{displaymath}
  It remains to sum the geometric series to obtain the result.
\end{proof}
We now want to obtain estimates for $\psi_\rho$ and $u_\rho$
restricted to $X$. We have the following result.
\begin{corollary}
  \label{cor:regX}
  Let us assume the regularity hypotheses of the previous proposition.
  Then we find that
  \begin{equation}
    \label{eq:reg}
    |\rho| \|\psi_\rho\|_{H^{\frac n2+k+\eps}(X)}
    + \|\psi_\rho\|_{H^{\frac n2+k+1+\eps}(X)} 
   \leq C  \|q\|_{H^{\frac n2+k+\eps}(X)}.
  \end{equation}
\end{corollary}
\begin{proof}
  On the bounded domain $X$, $\aver{x}$ is bounded above and below by
  positive constants. Since $q$ is compactly supported on $\Rm^n$, we
  obtain thanks to \eqref{eq:bdext} that
  \begin{displaymath}
     |\rho| \|\psi_\rho\|_{H^{\frac n2+k+\eps}(X)}
   \leq C \|q\|_{H^{\frac n2+k+\eps}(\Rm^n)}
    \leq C(X)  \|q\|_{H^{\frac n2+k+\eps}(X)}.
  \end{displaymath}
  Now we have 
  \begin{displaymath}
    \Delta \psi_\rho = -2\rho\cdot\nabla \psi_\rho -q(1+\psi_\rho).
  \end{displaymath}
  By elliptic regularity, with $X'$ a smooth domain in $\Rm^n$ such
  that $\bar X\subset X'$, we find for all $s=\frac n2+k+\eps$, that
  \begin{displaymath}
    \|\psi_\rho\|_{H^{s+1}(X)} \leq C 
        \|2\rho\cdot\nabla \psi_\rho -q(1+\psi_\rho)\|_{H^{s-1}(X')} +
   \|\psi_\rho\|_{H^s(X')}.
  \end{displaymath}
  The latter is bounded by
  $|\rho|\|\psi_\rho\|_{H^{s}(X')}+\|q\|_{H^{s-1}(X')}
  \|\psi_\rho\|_{H^{s-1}(X')}$
  since $s-1>\frac n2$ so that $H^{s-1}(X')$ is a Banach algebra. By
  using the above bound on $\|\psi_\rho\|_{H^{s}(X')}$, with $C(X)$
  replaced by the larger $C(X')$, we get the result.
\end{proof}
By Sobolev embedding, we have just proved the:
\begin{proposition}\label{prop:bdck}
  Under the hypotheses of Corollary \ref{cor:regX}, the restriction
  to $X$ of the CGO solution verifies that
  \begin{equation}
    \label{eq:boundCk}
    |\rho|\|\psi_\rho\|_{C^k(\bar X)} + \|\psi_\rho\|_{C^{k+1}(\bar X)}
    \leq C \|q\|_{H^{\frac n2+k+\eps}(X)}.
  \end{equation}
\end{proposition}
We recall that $q$ satisfies the constraint \eqref{eq:constq}.

We are now in a position to prove the main result of this section.
\begin{theorem}
  \label{thm:field}
  Let $u_{\rho_j}$ for $j=1,2$ be CGO solutions with $q$ such that
  \eqref{eq:constq} holds for both $\rho_j$ and $k\geq1$ and with
  $c_0^{-1}|\rho_1|\leq |\rho_2|\leq c_0|\rho_1|$ for some $c_0>0$.
  Then we have
  \begin{equation}
    \label{eq:drift}
    \hat\beta:= \dfrac{1}{2|\rho_1|} e^{-(\rho_1+\rho_2)\cdot x}
   \Big( u_{\rho_1}\nabla u_{\rho_2} - u_{\rho_2}\nabla u_{\rho_1}\Big)
    = \dfrac{\rho_1-\rho_2}{2|\rho_1|}
     + \hat h,
  \end{equation}
  where the vector field $\hat h$ satisfies the constraint
  \begin{equation}
    \label{eq:consthrho}
    \|\hat h\|_{C^k(\bar X)} \leq \dfrac{C_0}{|\rho_1|},
  \end{equation}
  for some constant $C_0$ independent of $\rho_{1,2}$.
\end{theorem}
\begin{proof}
  Some algebra shows that
  \begin{equation}
    \label{eq:h}
    \hat h =\dfrac{(\rho_1-\rho_2)}{2|\rho_1|}
   (\psi_{\rho_1}+\psi_{\rho_2}+\psi_{\rho_1}\psi_{\rho_2})
     +\dfrac{\nabla\psi_{\rho_2}(1+\psi_{\rho_1}) 
     -\nabla\psi_{\rho_1}(1+\psi_{\rho_2})}{2|\rho_1|}.
  \end{equation}
  We know from Proposition \ref{prop:bdck} that
  $|\rho_j||\psi_{\rho_j}|$ and $|\nabla\psi_{\rho_j}|$ are bounded in
  $C^k(\bar X)$ for $j=1,2$. This concludes the proof of the theorem.
\end{proof}

%
%
\subsection{Construction of vector fields and uniqueness result}
\label{sec:vector fields}

Let us consider two internal complex-valued data $d_{1,2}(x)$ obtained
as follows.  We assume that we can impose the complex-valued boundary
conditions $g_{1,2}\in C^{k,\alpha}(\pa X;\Cm)$ and define the
solution $u_{1,2}$ of
\begin{equation}
  \label{eq:u}
  \Delta u_j + q u_j =0, \quad X,\qquad
   u_j = g_j, \quad \pa X, \qquad j=1,2.
\end{equation}
Note that the real and imaginary parts of $u_{1,2}$ may be solved
independently since \eqref{eq:u} is a linear equation.  We then assume
that we have access to the complex-valued internal data $d_j= \mu u_j$
on $X$ for $j=1,2$, where $u_{1,2}$ are the solutions of \eqref{eq:u}
with boundary conditions $g_{1,2}$. We recall that $\mu\in
C^{k+1}(\bar X)$ and is bounded above an below by positive constants.
We verify that
\begin{displaymath}
  u_1\Delta u_2 - u_2 \Delta u_1 =0.
\end{displaymath}
Introducing $\nu=\frac1{\mu}$, which is well defined since $\mu$ is
bounded away from $0$, and using $u_j=\nu d_j$, we obtain that
\begin{displaymath}
  2(d_1\nabla d_2-d_2\nabla d_1)\cdot\nabla\nu
  +(d_1\Delta d_2-d_2\Delta d_1)\nu=0.
\end{displaymath}
This is equivalent to
\begin{equation}
  \label{eq:formu}
  \check\beta_d\cdot\nabla\mu + \check\gamma_d \mu =0,
\end{equation}
where 
\begin{equation}
  \label{eq:betagammad}
   \begin{array}{rcl}
  \check\beta_d &:=&  \chi(x) (d_1\nabla d_2-d_2\nabla d_1) \\
  \check\gamma_d &:=&  \dfrac12\chi(x) (d_2\Delta d_1-d_1\Delta d_2) = 
   \dfrac{-\check\beta_d\cdot\nabla\mu}{\mu}.
   \end{array}
\end{equation}
Here, $\chi(x)$ is a smooth known complex-valued function with
$|\chi(x)|$ uniformly bounded from below by a positive constant on
$\bar X$.  Note that by assumption on $\mu$, we have that
$\check\beta_d\in (C^k(\bar X;\Cm))^n$ and $\check\gamma_d\in C^k(\bar
X;\Cm)$.

A methodology for the reconstruction of $(\mu,q)$ is therefore as
follows: we first reconstruct $\mu$ using the real part {\em or} the
imaginary part of \eqref{eq:formu} for then $\Re\check\beta_d$ and
$\Im\check\beta_d$ are real-valued vector fields since $\mu=d/g$ is
known on $\pa X$. When $\mu$ is reconstructed, this gives us explicit
reconstructions for $u_{1,2}=d_{1,2}/\mu$ and we may then reconstruct
$q$ from the Schr\"odinger equation. Such a method provides a unique
reconstruction provided that the integral curves of (the real part or
the imaginary part of) $\check\beta_d$ join any point in $x$ to a
point $x_0(x)\in\pa X$, where $\mu$ is known. We thus need the vector
field $\check\beta$ to satisfy such properties.  CGO solutions will
allow us to construct families of vector fields $\check\beta_d$ with
the required properties.

Let us consider two CGOs $u_{\rho_{1,2}}$with parameters $\rho_{1,2}$.
Let $d_{1,2}$ be the complex-valued corresponding internal data.  Let
us decompose as before
\begin{displaymath}
  u_{\rho_j}(x) = e^{\rho_j\cdot x} (1+\psi_{\rho_j}(x)),\qquad
  \nabla u_{\rho_j}(x) = e^{\rho_j\cdot x}
  \big((1+\psi_{\rho_j})\rho_j+\nabla\psi_{\rho_j}\big).
\end{displaymath}
Let us choose $\chi(x)=e^{-(\rho_1+\rho_2)\cdot x}$ in
\eqref{eq:betagammad}.  Then we find after some algebra that
$\check\beta_d$ in \eqref{eq:betagammad} is given by
\begin{equation}
  \label{eq:betadrho12}
  \check\beta_d = \mu^2
   \Big((\rho_1-\rho_2)(1+\psi_{\rho_1})(1+\psi_{\rho_2})
   +\nabla\psi_{\rho_2}(1+\psi_{\rho_1})
    -\nabla\psi_{\rho_1}(1+\psi_{\rho_2})\Big).
\end{equation}
We may then define
\begin{equation}
  \label{eq:checkbeta}
  \check\beta := \dfrac{1}{2|\rho_1|} \check \beta_d = 
  \mu^2 \dfrac{\rho_1-\rho_2}{2|\rho_1|} + \mu^2 \hat h, \qquad
  \check\gamma :=  \dfrac{1}{2|\rho_1|} \check \gamma_d,
\end{equation}
where $\hat h$ is defined as in \eqref{eq:h}. Then, we deduce from
Theorem \ref{thm:field} that $|\rho_1|\mu^2\hat h$ is bounded
uniformly in $C^k(\bar X;\Cm)$.  When $|\rho_1|$ is sufficiently
large, then $\check\beta$ is close to
$\mu^2\frac{\rho_1-\rho_2}{2|\rho_1|}$, which is a non-vanishing
vector when $\rho_1\not=\rho_2$.  Provided that the real part or the
imaginary part of $\check\beta$ does not vanish, then \eqref{eq:formu}
gives an equation for $\mu$ that can be uniquely solved since $\mu$ is
known on $\pa X$.

\medskip

Note that the data $d_k$ are complex valued.  The only possibility to
construct two different complex valued data with two real valued data
is to assume that $d_2=\bar{d_1}$, the complex conjugate of $d_1$.
For the construction of CGOs, this implies that we choose
$\rho_2=\overline{\rho_1}$. Indeed, we verify that
$\overline{u_\rho}=u_{\bar\rho}$ since for $\rho=\k+i\k^\perp$ with
$|\k|=|\k^\perp|$ and $\k\cdot\k^\perp=0$, we have
$\overline{e^{\rho\cdot x}}=e^{\bar \rho\cdot x}$ and
$\overline{\psi_\rho}=\psi_{\bar\rho}$ by uniqueness of the solution
to the equation satisfied by $\psi_\rho$. This implies then that
$\check\beta$ defined in \eqref{eq:checkbeta} with
$\rho_2=\bar{\rho_1}$ is given by
\begin{equation}
  \label{eq:betabar}
  \check\beta_\rho = i\mu^2 \k^\perp + \mu^2\hat h.
\end{equation}
As soon as $|\rho|>C_0$ so that $\|\hat h\|_{C^0(\bar X)}<1$, we
obtain that any point in $X$ is connected to a point in $\pa X$ by an
integral curve of $\beta_\rho:=\Im\check\beta_\rho$.

\medskip

Note that $u_{\rho}$ solves \eqref{eq:u} with the unknown boundary
condition $u_{\rho|\pa X}\in C^{k,\alpha}(\pa X;\Cm)$ for some
$\alpha>\frac12$ since $u_\rho$ is known to be a little more regular
than being of class $C^{k+1}(\bar X;\Cm)$ by construction (since $\eps>0$).

Let us now define boundary conditions $g\in C^{k,\alpha}(\pa X;\Cm)$
such that 
\begin{equation}
  \label{eq:closeg}
  \|g-u_{\rho|\pa X}\|_{C^{k,\alpha}(\pa X;\Cm)} \leq \epsilon,
\end{equation}
for some $\epsilon>0$ sufficiently small.
Let $u$ be the solution of \eqref{eq:u} with $g$ as in
\eqref{eq:closeg}. By elliptic regularity, we thus have
\begin{equation}
  \label{eq:errorurho}
  \|u-u_\rho\|_{C^{k+1}(\bar X;\Cm)} \leq C\epsilon,
\end{equation}
for some positive constant $C$. Define the complex valued internal
data $d=\mu u$. Since $\mu\in C^{k+1}(\bar X)$, we deduce that
\begin{equation} \label{eq:approxd}
  \|d-d_\rho\|_{C^{k+1}(\bar X;\Cm)} \leq C_0\epsilon,
\end{equation}
for $C_0>0$. Once $d$ is constructed, define $d_1=d$ and $d_2=\bar d$
and define $\check\beta_d$ and $\check\mu_d$ as in
\eqref{eq:betagammad} with $\chi(x)=e^{-2\k\cdot x}$ and the
normalized quantities $\check\beta$ and $\check\gamma$ as in
\eqref{eq:checkbeta}. Note that $\chi(x)$ is positive and bounded on
$\bar X$.

Let us define 
\begin{equation}
  \label{eq:beta}
  \beta := \Im\check\beta = \dfrac{1}{2|\k|} \Im \check\beta_d,\qquad
  \gamma := \Im\check\gamma = \dfrac{1}{2|\k|} \Im \check\gamma_d.
\end{equation}
Thanks to \eqref{eq:approxd} and \eqref{eq:betabar}, we obtain the
error estimate
\begin{equation}
  \label{eq:controlbeta}
  \|\beta-\mu^2\hat\k\|_{C^{k}(\bar X)} \leq C\dfrac{1+\epsilon}{|\k|}.
\end{equation}
As a consequence, as soon as $|\k|$ is sufficiently large and
$\epsilon$ sufficiently small, we obtain that
$\beta\cdot\hat\k\geq\zeta>0$ so that any point $x\in X$ is mapped to
a point in $\pa X$ in a time less than $|\zeta|^{-1}{\rm diam}(X)$ by
an integral curve of $\beta$.

Moreover, we have the equation with real-valued coefficients:
\begin{equation}
  \label{eq:formu1}
  \beta\cdot\nabla\mu + \gamma \mu =0.
\end{equation}
Since $\mu=d/g$ is known on $\pa X$, this equation provides a unique
reconstruction for $\mu$.

\medskip

Let us define the set of parameters 
\begin{equation}
  \label{eq:P}
   \begin{array}{rcl}
  \mathcal P&\!\!  =\!\!
  & \Big\{(\mu,q)\in C^{k+1}(\bar X)\times H^{\frac n2+k+\eps}(X);
   \,\, 0\mbox{ not an eigenvalue of }\Delta +q, \\
  &&\qquad
   \|\mu\|_{C^{k+1}(\bar X)} + \|q\|_{H^{\frac n2+k+\eps}(X)} \leq P
  <\infty  \Big\}.
  \end{array}
\end{equation}
The above construction of the vector field allows us to obtain the 
following uniqueness result.
\begin{theorem}
  \label{thm:uniq}
  Let $X$ be a bounded, open subset of $\Rm^n$ with boundary $\pa X$
  of class $C^2$.  Let $(\mu,q)$ and $(\tilde\mu,\tilde q)$ be two
  elements in $\mathcal P$.  Let $\k\in\Rm^n$ with $|\k|\geq|\k_0|$
  and $|\k_0|$ sufficiently large and define $\rho=\k+i\k^\perp$ so
  that $\rho\cdot\rho=0$. Let $u_\rho$ be the corresponding $CGO$ for
  $q$ and $u$ constructed as above with $d=\mu u$ and with $\epsilon$
  sufficiently small. Let $\tilde d$ be constructed similarly with the
  parameters $(\tilde\mu,\tilde q)$.
  
  Then $d=\tilde d$ implies that $(\mu,q)=(\tilde\mu,\tilde q)$.
\end{theorem}
\begin{proof}
  Since the two measurements $d=\tilde d$, we have that $\mu$ and
  $\tilde\mu$ solve the same equation \eqref{eq:formu1}. Since
  $\mu=\tilde\mu=d/g$ on $\pa X$, we deduce that $\mu=\tilde\mu$ since
  the integral curves of $\beta$ map any point $x\in X$ to the
  boundary $\pa X$. More precisely, consider the flow $\varphi_x(t)$
  associated to $\beta$, i.e., the solution to
  \begin{equation}
  \label{eq:ODE}
  \dot\varphi_x(t) = \beta(\varphi_x(t)), \qquad
   \varphi_x(0)=x \in \bar X.
  \end{equation}
  By the Picard-Lindel\"of theorem, the above equations admit unique
  solutions since $\beta$ is of class $C^1$.  And by hypothesis on
  $\beta$ since $\beta\cdot\hat\k\geq\zeta>0$ for $|\k|$ sufficiently
  large, any point $x$ is mapped to two points (for positive and
  negative values of $t$) on $\pa X$ by the flow $\varphi_x(t)$ in a
  time less than $\zeta^{-1}{\rm diam}(X)$.  For $x\in X$, let us
  define $x_\pm(x)\in\pa X$ and $\pm t_\pm(x)>0$ such that
  \begin{equation} \label{eq:xtpm}
    \varphi_x(t_\pm(x)) = x_\pm(x) \in \pa X.
  \end{equation}
  Then by the method of characteristics, $\mu(x)$ solution of
  \eqref{eq:formu1} is given by 
  \begin{equation}
  \label{eq:solode}
  \mu(x) = \mu_0(x_\pm(x))e^{-\int_0^{t_\pm(x)}\gamma(\varphi_x(s))ds}\,.
  \end{equation}
  The solution $\tilde\mu(x)$ is given by the same formula since
  $\varphi_x(t)=\tilde\varphi_x(t)$ so that $\tilde\mu=\mu$.  This
  implies that $u=\tilde u$ since $d=\tilde d$. It remains to use the
  equation for $u$ to deduce that $q=\tilde q$ on the domain where
  $u\not=0$. By unique continuation, $u$ cannot vanish on an open set
  in $X$ different from the empty set for otherwise $u$ vanishes
  everywhere and this is impossible to satisfy the boundary
  conditions. This shows that the set $F\subset X$ where $|u|>0$ is
  open and $\bar F=\bar X$ since the complement of $\bar F$ has to be
  empty. By continuity, this shows that $q$ is known on $\bar X$.
\end{proof}

The above result shows that there exists an open set of boundary
conditions $g$ close to $u_{\rho|\pa X}$ so that data $d_1=d$ and
$d_2=\bar d$ obtained from one complex-valued solution $u$ or
equivalently from two real valued solutions $\Re u$ and $\Im u$,
uniquely determine the parameters $(\mu,q)$. A more explicit
characterization of the open set of illuminations is lacking.
However, we observe that larger values of $q$ require larger values of
$|\k|$ in order to straighten the vector field $\beta$.  Although
\eqref{eq:closeg} seems to be independent of $\k$ and $\rho$, in fact
$u_\rho$ itself grows exponentially with $|\k|=|\rho|$ so that $g$ has
to be in the $\epsilon$ vicinity of an exponentially growing quantity.
This means that $|\k|$ has to be sufficiently large that the field
$\beta$ is sufficiently flat while at the same time not so large that
the imposed illuminations become physically infeasible.
 
\medskip

The above uniqueness result does not guaranty stability in the
reconstruction. We easily verify that the construction provides
stability of the reconstruction of $\mu$ in most of the domain $X$.
However, small changes in the data may generate small changes in the
field $\beta$.  This in turn may significantly modify the value of the
reconstructed function $\mu$ at points where $\beta$ is ``almost''
tangent to the boundary $\pa X$. We will see below that under some
geometric constraints of sufficient convexity of $X$, the above
procedure provides a stable reconstruction of the parameters
$(\mu,q)$. When such conditions are not met, we can still obtain
stability by acquiring more measurements. Indeed, if a sufficient
number of vector fields $\beta$ can be constructed at every point so
that the span of these vector fields is exactly $\Rm^n$, then we face
a significantly more favorable situation. We now consider such a case
where $2n$ real-valued measurements are available.  Later, we will
derive stability results in the two-measurement setting under
additional geometric constraints.

%
%
%
%
%
\subsection{ISID with $2n$ real-valued internal data}
\label{sec:multmeas}

Let us consider first the setting in which we can access $2n$
real-valued internal data viewed as $n$ complex-valued internal data
(since the measurements are linear in $u$, we can measure the real and
imaginary parts separately).

Let us define $\k_j=|\k|e_j$ where $(e_1,\ldots,e_n)$ is an
orthonormal basis. We define the complex vectors 
\begin{equation}
  \label{eq:rhoj}
  \rho_j = \k_1+i\k_j, \quad 2\leq j\neq n,\qquad
  \rho_1=  -\k_1-i\k_2 = -\rho_2.
\end{equation}
Let $u_{\rho_j}$ be the corresponding CGOs. We choose boundary conditions
$g_j$ such that
\begin{equation}
  \label{eq:closegj}
  \|g_j-u_{\rho_j|\pa X}\|_{C^{k,\alpha}(\pa X;\Cm)} \leq \epsilon,
\end{equation}
for $\epsilon$ sufficiently small. We define $u_j$ as the solutions to
\eqref{eq:u} with boundary conditions $g_j$. These are $n$
complex-valued solutions whose real and imaginary parts consist of
$2n$ real-valued solutions. For $1\leq j\leq n$, we define $d_j=\mu
u_j$. We now construct the $n$ vector field $\beta_j$. For $2\leq
j\leq n$, the real-valued vector fields and scalar terms are
constructed as in the preceding section; for $j=1$, the vector field
is constructed by using $\rho_2=-\rho_1$:
\begin{equation}
  \label{eq:betagammaj}
  \begin{array}{rclrcl}
       \beta_1 &=&\dfrac{1}{2|\k|} \Re\big(d_2\nabla d_1-d_1\nabla d_2\big),
   \,\,&  \gamma_1 &=& \dfrac{1}{4|\k|} \Re
   \big(d_1\Delta d_2-d_2\Delta d_1\big), \\
    \beta_j &=& \dfrac{e^{-2\k_1\cdot x}}{2|\k|} \Im\big(d_j\nabla \bar{d_j}
    -\bar {d_j}\nabla d_j\big),\,\,&
    \gamma_j &=&  \dfrac{e^{-2\k_1\cdot x}}{4|\k|} \Im
       \big(\bar{d_j}\Delta d_j-d_j\Delta{\bar d_j}\big),
  \end{array}
\end{equation}
for $2\leq j\leq n$. As in the preceding section, we verify that
\begin{equation}
  \label{eq:errorbetas}
  \|\beta_j-\mu^2\hat\k_j\|_{C^{k}(\bar X)} \leq C\dfrac{1+\epsilon}{|\k|}.
\end{equation}

For $|\k|$ sufficiently large, and thanks to the bound $\mu_0^{-1}\leq
\mu\leq \mu_0$, we obtain that at each point $x\in X$, the vectors
$\beta_{j}(x)$ form a basis. Moreover, the matrix $a_{ij}$ such that
$\beta_{j}= \sum a_{jk} e_k$ is an invertible matrix with inverse of
class $C^k(\bar X)$. In other words, we have constructed a
vector-valued function $\Gamma(x)\in (C^k(\bar X))^n$ such that
\eqref{eq:formu1} may be recast as
\begin{equation}
  \label{eq:gradnu}
  \nabla \mu + \Gamma(x) \mu =0.
\end{equation}
Finally, the construction of $\Gamma$ is stable under small perturbations
in the data $d_j$. Indeed, invertibility of $a_{jk}$ is ensured for
vector fields close to $\beta_{j}$. Let $\Gamma$ and $\tilde\Gamma$
be two vector fields constructed from knowledge of two sets of
internal data $d=\{d_{j},1\leq j\leq n\}$ and 
$\tilde d=\{\tilde d_{j},1\leq j\leq n\}$. Then we find that
\begin{equation}
  \label{eq:stabGamma}
  \|\Gamma-\tilde\Gamma\|_{(C^{k}(\bar X))^n}
   \leq C \|d-\tilde d\|_{(C^{k+1}(\bar X;\Cm))^{n}},
\end{equation}
provided the right-hand side is sufficiently small.

Let us now assume that $X$ is connected (otherwise, the method applies
to each connected component) and $\mu$ is known and equal to
$\mu_0=d/g$ for some point $x_0\in \pa X$.  In other words, we want to
solve the over-determined problem
\begin{equation}
  \label{eq:nuoverdet}
  \nabla \mu + \Gamma(x) \mu =0 ,\quad X,\qquad
   \mu(x_0) = \mu_0(x_0),\quad x_0 \in\pa X.
\end{equation}
Let $x\in X$ be an arbitrary point and assume that $X$ is bounded and
connected and $\pa X$ is smooth.  Then we find a smooth curve that
links $x$ to the point $x_0\in\pa X$.  Restricted to this curve,
\eqref{eq:nuoverdet} becomes a stable ordinary differential equation.
The solution of the ordinary differential equation is then stable with
respect to modifications in $\Gamma$ (the curve between $x$ and $x_0$
is kept constant).  The solution $\mu$ then clearly inherits the
smoothness of $\Gamma(x)$ directly from \eqref{eq:nuoverdet}.
Moreover, since $\mu_0(x_0)-\tilde\mu_0(x_0)$ (with
$\tilde\mu_0=\tilde d/g$ on $\pa X$) is small and equation
\eqref{eq:stabGamma} is stable with respect to changes in the value of
$\mu_0(x_0)$, we deduce that the reconstruction of $\mu$ is stable
with respect to perturbations in $d$.

We may thus state the main result of this section:
\begin{theorem}
  \label{thm:stab2n}
  Let $k\geq1$.  We assume that we have access to $n$ well-chosen
  complex-valued measurements and that $(\mu,q)$ and
  $(\tilde\mu,\tilde q)$ are elements in $\mathcal P$. Under the
  hypotheses outlined above, and provided that $\|u_{\rho_j|\pa
    X}-g_j\|_{C^{k,\alpha}(\bar X;\Cm)}$ is sufficiently small, then we
  have the following stability result:
  \begin{equation}
  \label{eq:stabmuq}
   \|\mu-\tilde \mu\|_{C^{k}(\bar X)} 
     +\|q-\tilde q\|_{C^{k-2}(\bar X)} 
     \leq C \|d-\tilde d\|_{(C^{k+1}(\bar X))^{2n}},
  \end{equation}
\end{theorem}
\begin{proof}
  The inequality for $\mu-\tilde\mu$ is a direct consequence of the
  results proved above. This provides a stability result for
  $\nu=\mu^{-1}$ and for $u_{j}=\nu d_j$ from the data $d_{j}$. We
  thus have a stability result for $\Delta u_j=-u_{j} q$ and hence the
  above stability result for $u_{j}(q-\tilde q)$ since $(u_{j}-\tilde
  u_{j})\tilde q$ is small. 
  
  Now, $u_\rho=e^{\rho\cdot x}(1+\psi_\rho)$ does not vanish on $X$
  when $|\rho|$ is sufficiently large since $|\rho|\psi_\rho$ is
  bounded. When the boundary condition $g_{j}-u_{\rho_j|\pa X}$ is
  small, then by the maximum principle, $u_j$ does not vanish on $X$
  either.  This means that either its real part or its imaginary part
  does not vanish everywhere in $X$. This provides control of
  $q-\tilde q$ in $X$ as given in \eqref{eq:stabmuq}.
\end{proof}

%
%
%
%
%
\subsection{Vector fields and stability of solutions}
\label{sec:stabvf}

The above construction allows one to stably reconstruct the two
functions $\mu$ and $q$ provided that we have constructed $J=2n$
well-chosen real-valued boundary conditions and collected $2n$
corresponding internal data. We now return to the reconstruction of
$\mu$ and $q$ in the presence of $J=2$ well-chosen real-valued
internal data.  Such internal data are obtained as Theorem
\ref{thm:uniq}.

We recall that $\k$ is fixed and $\rho=\k+i\k^\perp$. We define $u$ as
the solution to \eqref{eq:u} with $g$ close to $u_{\rho|\pa X}$.  The
complex-valued internal data are then $d=\mu u$.  The vector field
$\beta$ and the scalar $\gamma$ are then given by
\begin{equation}
  \label{eq:betagamma0}
  \begin{array}{rclrcl}
    \beta &=& \dfrac{e^{-2\k\cdot x}}{2|\k|} \Im\big(d\nabla \bar{d}
    -\bar {d}\nabla d\big),\,\,&
    \gamma &=&  \dfrac{e^{-2\k\cdot x}}{4|\k|} \Im
       \big(\bar{d}\Delta d-d\Delta{\bar d}\big),
  \end{array}
\end{equation}
and we verify that
\begin{equation}\label{eq:bgmu}
  \beta\cdot\nabla\mu +\gamma\mu =0 \qquad X.
\end{equation}
As earlier, we verify that
\begin{equation}
  \label{eq:errorbeta0}
  \|\beta-\mu^2\hat\k^\perp\|_{C^{k}(\bar X)} \leq C\dfrac{1+\epsilon}{|\k|}.
\end{equation}
As a consequence, the integral curves of $\beta$ given by
$\varphi_x(t)$ map any point $x\in X$ to two points on $\pa X$ when
$|\k|$ is sufficiently large as was mentioned earlier.

However, the stability of equation \eqref{eq:bgmu} with respect to
changes in $\beta$ and $\gamma$ is not as good as in the presence of
$n$ complex internal data.  The stability of the reconstruction
degrades for points $x$ close to $x_0\in\pa X$ where
$\beta_\mu(x_0)\cdot n(x_0)$ is close to $0$. The stability of the
reconstruction of $\mu$ will however be good when $X$ is a convex
domain with ``sufficient'' convexity as established in Hypothesis
\ref{hyp:R0}.  We prove the following result:

\begin{proposition}
  \label{prop:regmu}
  Let $k\geq1$.  Let $\mu$ and $\tilde\mu$ be solutions of
  \eqref{eq:bgmu} corresponding to coefficients $(\beta,\gamma)$ and
  $(\tilde\beta,\tilde\gamma)$, respectively, where
  \begin{displaymath}
    \beta= \mu^2 \hat \k^\perp + \dfrac{1}{|\k|} h,\qquad
    \tilde\beta= \mu^2\hat \k^\perp + \dfrac{1}{|\k|} \tilde h,
  \end{displaymath}
  for $h$, $\gamma$, $\tilde h$, and $\tilde\gamma$ bounded in
  $C^{k}(\bar X)$.
  
  Let us assume that $\mu_{|\pa X}=\mu_0$ and $\tilde\mu_{|\pa
    X}=\tilde\mu_0$ on $\pa X$ for some functions
  $\mu_0,\tilde\mu_0\in C^{k}(\pa X)$. Let us assume that $X$ is
  sufficiently convex so that Hypothesis \ref{hyp:R0} holds for some
  $R<\infty$. We also assume that $|\k|\geq \k_0$ is sufficiently
  large.  Then there is a constant $C$ such that
  \begin{equation}
    \label{eq:stabnu2}
    \begin{array}{rcl}
    \|\mu-\tilde\mu\|_{C^{k-1}(\bar X)} &\leq& C
    \|\mu_0\|_{C^{k}(\pa X)} 
    \big(\|\beta-\tilde\beta\|_{C^{k-1}(\bar X)}
    + \|\gamma-\tilde\gamma\|_{C^{k-1}(\bar X)}\big)\!\!\!\! \\[3mm]
    & +& C \|\mu_0-\tilde\mu_0\|_{C^{k}(\pa X)}.
    \end{array}
  \end{equation}
\end{proposition}

\begin{proofof} {\it Proposition \ref{prop:regmu}}.
  We recall that $\varphi_x(t)$ is the flow defined in \eqref{eq:ODE}
  and that $x_\pm(x)$ and $t_\pm(x)$ are defined in \eqref{eq:xtpm}.
  By the method of characteristics, $\mu(x)$ solution of
  \eqref{eq:bgmu} is given by
\begin{equation}
  \label{eq:solode2}
  \mu(x) = \mu_0(x_\pm(x))e^{-\int_0^{t_\pm(x)}\gamma(\varphi_x(s))ds}\,.
\end{equation}
The solution $\tilde\mu(x)$ is given similarly. 
We first assume that $\tilde\mu_0=\mu_0$.

From the equality
\begin{displaymath}
  \varphi_x(t)-\tilde\varphi_x(t) = \dint_0^t
   [\beta(\varphi_x(s))-\tilde\beta(\tilde\varphi_x(s))] ds,
\end{displaymath}
and using the Lipschitz continuity of $\beta$ and Gronwall's lemma, we
thus deduce the existence of a constant $C$ such that
\begin{displaymath}
  |\varphi_x(t)-\tilde\varphi_x(t)| \leq Ct \|\beta-\tilde\beta\|_{C^0(X)}
\end{displaymath}
uniformly in $t$ knowing that all characteristics exit $X$ in finite
time and provided that $\varphi_x(t)$ and $\tilde\varphi_x(t)$ are in
$\bar X$.

Such estimates are stable with respect to modifications in the initial
conditions. Let us define $W(t)=D_x\varphi_x(t)$. Then classically,
$W$ solves the equation $\dot W= D_x\beta(\varphi_x)W$ with $W(0)=I$
and by using Gronwall's lemma once more, we deduce that
\begin{displaymath}
  |W-\tilde W|(t) \leq Ct \|D_x\beta-D_x\tilde\beta\|_{C^0(\bar X)},
\end{displaymath}
for all times provided that $\varphi_x(t)$ and $\tilde\varphi_x(t)$
are in $\bar X$. As a consequence, since $\beta$ and $\tilde\beta$
are of class $C^{k}(\bar X)$, then we obtain similarly that:
\begin{displaymath}
  |D_x^{k-1}\varphi_x(t)-D_x^{k-1}\tilde\varphi_x(t)| 
   \leq Ct \|\beta-\tilde\beta\|_{C^{k-1}(X)},
\end{displaymath}
and this again for all times provided that $\varphi_x(t)$ and
$\tilde\varphi_x(t)$ are in $\bar X$.

However, this does not imply that $x_+(x)$ is close to $\tilde
x_+(x)$.  When $\pa X$ is flat for instance, we may very well have
that $x_+(x)$ is such that $n(x_+(x))\cdot \dot\varphi_x(t_+(x))$ is
very small and that $x_+(x)-\tilde x_+(x)$ is arbitrarily large if
$\tilde\beta$ is parallel to the surface $\pa X$ for instance.  This
behavior, however, cannot occur when both $\beta$ and $\tilde\beta$
are sufficiently flat, which is the case when $|\k|$ is sufficiently
large, and when $\pa X$ is sufficiently curved, which is obtained from
the existence of $R<\infty$ in Hypothesis \ref{hyp:R0}. In such a
setting, we can obtain the following result:
\begin{lemma}
  \label{lem:bds}
  Let $k\geq1$ and assume that $\beta$ and $\tilde\beta$ are
  $C^{k}(\bar X)$ vector fields that are sufficiently flat, i.e., that
  $|\k|$ is sufficiently large. Let us assume that $\pa X$ is
  sufficiently convex so that Hypothesis \ref{hyp:R0} holds for some
  $R<\infty$. Then we have that
  \begin{equation}
    \label{eq:boundxt}
    \|x_+-\tilde x_+\|_{C^{k-1}(\bar X)}+\|t_+-\tilde t_+\|_{C^{k-1}(\bar X)}
   \leq C \|\beta-\tilde\beta\|_{C^{k-1}(\bar X)},
  \end{equation}
  where $C$ is a constant that depends on $|\k|$ and $R$.
\end{lemma}

The above lemma is mostly a consequence of the following result:
\begin{lemma}
  \label{lem:propxt}
  Let $C_0$ be the constant defined such that
  \begin{displaymath}
    |\ddot\varphi_x(t)| = |\nabla\beta(\varphi_x(t)) \beta (\varphi_x(t))|
   \leq \dfrac{C_0}{|\k|}. 
  \end{displaymath}
  Let $t_M$ be the maximal time spent by any trajectory in $X$, which
  we know is bounded.  Assume that $|\k|$ is sufficiently large that
  for all $x\in X$,
  \begin{displaymath}
    \Big(\dfrac{C_0R}{|\k|} +\dfrac{C_0^2t_M^2}{4|\k|^2}\Big) 
   \dfrac{1}{|\beta(x_+(x))|^2} = \rho<1.
  \end{displaymath}
  Then we have that for all $x\in X$,
  \begin{equation}
    \label{eq:boundt}
    t_+(x)\leq\frac{2R}{|\beta(x_+(x))|^2(1-\rho)}n(x_+(x))\cdot\beta(x_+(x)).
  \end{equation}
  In other words, the vector field $\beta(x_+(x))$ is close to being
  tangent to $\pa X$ only when the time spent in $X$ is small.
  
  Let now $x_0\in X$ and $x\in \pa X$ and define $v_0=\beta(x_0)$.
  Assume moreover that
  \begin{equation}
    \label{eq:bds}
    |x-x_0|\leq C_1 \delta^2,\qquad
    |n(x)\cdot v_0| \leq C_2 \delta.
  \end{equation}
  Then we have
  \begin{equation}
    \label{eq:bdtp}
    t_+(x_0) \leq C_3\delta,
  \end{equation}
  for some constant $C_3>0$ independent of $x_0$.  In other words, a
  trajectory close to $\pa X$ and almost tangent to $\pa X$ exits $X$
  in a short time.
\end{lemma}

We postpone the proof of the above two lemmas to the end of the
section.  Let us conclude the proof of the proposition. We recall that
\begin{displaymath}
  \mu(x) = \mu_0(x_+(x)) e^{-\int_0^{t_+(x)}\gamma(\varphi_x(s))ds},
\end{displaymath}
with a similar expression for $\tilde\mu$. Since $x\to e^{-x}$ is
smooth, by the Leibniz rule it is sufficient to prove the stability
result for $\mu_0(x_+(x))$ and for
$\int_0^{t_+(x)}\gamma(\varphi_x(s))ds$. It is clear from the above
lemmas that
\begin{displaymath}
  \begin{array}{l}
  \|\mu_0(x_+(x))-\mu_0(\tilde x_+(x)) \|_{C^{k-1}(\bar X)} 
  \leq \|\mu_0\|_{C^{k}(\pa X)} \|x_+-\tilde x_+\|_{C^{k-1}(\bar X)}
  \\ 
  \leq C \|\mu_0\|_{C^{k}(\pa X)} \|\beta-\tilde\beta\|_{C^{k-1}(\bar X)}.
  \end{array}
\end{displaymath}

Let us now assume without loss of generality that $\tilde t_+(x)\geq
t_+(x)$. Then we have
\begin{displaymath}
  \dint_0^{t_+(x)} \hspace{-.5cm} \big(\gamma(\varphi_x(s))-
   \tilde\gamma(\tilde\varphi_x(s))\big) ds
  =\dint_0^{t_+(x)} \hspace{-.5cm}
  \big(\gamma(\varphi_x(s))-\gamma(\tilde\varphi_x(s))
    + (\gamma-\tilde\gamma)(\varphi_x(s))\big)ds.
\end{displaymath}
We verify that the above expression has $k-1$ derivatives uniformly
bounded since (i) $x\to t_+(x)$ is $C^{k-1}(\bar X)$; (ii) $\gamma$ has
$C^{k}$ derivatives bounded on $\bar X$; (iii)
$(\varphi_x-\tilde\varphi_x)(s)$ has $k-1$ derivatives bounded by
$\|\beta-\tilde\beta\|_{C^{k-1}(\bar X)}$.

It thus remains to handle the term 
\begin{displaymath}
  \upsilon(x) = \int_{t_+(x)}^{\tilde t_+(x)}
    \tilde\gamma(\tilde\varphi_x(s)) ds.
\end{displaymath}
The function $x\to\tilde\gamma(\tilde\varphi_x(s))$ is of class
$C^{k-1}(\bar X)$ by regularity of the flow and because $\tilde\gamma$ is
of class $C^{k}(\bar X)$. Derivatives of order $k-1$ of $\upsilon(x)$
thus involve terms of size $\tilde t_+(x)- t(x)$ and terms of the form
\begin{displaymath}
  D_x^m \Big( \tilde t_+ D_x^{k-1-m} \tilde\gamma(\tilde\varphi_x(\tilde t_+))
   - t_+ D_x^{k-1-m}\tilde\gamma(\tilde\varphi_x(t_+))\Big),\quad
  0\leq m\leq k-1.
\end{displaymath}
Because $\tilde\beta$ is of class $C^{k}(\bar X)$, then so is
$x\to \tilde\gamma (\tilde\varphi_x(s))$. Since the latter function
has $k-1$ derivatives that are Lipschitz continuous, we thus find that
\begin{displaymath}
  |D_x^{k-1} \upsilon (x)| \leq C  \|\tilde t_+- t_+ \|_{C^{k-1}(\bar X)}.
\end{displaymath}
This concludes the proof of the proposition when $\tilde\mu_0=\mu_0$.
Applying Lemma \ref{lem:bds} as before, we verify that
\begin{displaymath}
  \big|D_x^{k-1}\big[(\tilde\mu_0(x_\pm(x)) - \mu_0(x_\pm(x)))
   e^{-\int_0^{t_\pm(x)}\gamma(\varphi_x(s))ds}\big]\big|
  \leq C \|\mu_0-\tilde\mu_0\|_{C^{k}(\bar X)}.
\end{displaymath}
By the triangle inequality, we deduce the error estimate on
$\mu-\tilde\mu$ described in the proposition.
\end{proofof}

\begin{proofof}{\it Lemma \ref{lem:bds}}.
  Let us assume without loss of generality that $t_+(x)\leq \tilde t_+(x)$.
  We have seen that 
  \begin{displaymath}
    |\varphi_x(t_+(x))-\tilde\varphi_x(t_+(x))| 
    \leq Ct_+(x) \delta,\qquad 
      \delta= \|\beta-\tilde\beta\|_{C^0(\bar X)}.
  \end{displaymath}
  We also have that
  $|\beta(\vp_x(t_+))-\tilde\beta(\tilde\vp_x(t_+))|\leq C_2\delta$.
  From Lemma \ref{lem:propxt}, we know that $\beta(x_+(x))\cdot
  n(x_+(x)) \geq 2C_1 t_+(x)$ for some constant $C_1>0$.
  
  Let us assume first that $C_1 t_+ \geq C_2\delta$ so that $v_0\cdot
  e_1 \geq C_1 t_+(x)$ where $x_0=\tilde\vp_x(t_+)$,
  $v_0=\tilde\beta(x_0)$, and $e_1=n(x_+(x))$. We want to show that
  the integral curve of $\tilde\beta$ starting at $x_0$ at time $0$
  with velocity $v_0$ exits $X$ in a time of order $\delta$ so that
  $x_+(x)-\tilde x_+(x)$ is of order $\delta$.
  
  In an appropriate system of coordinates, we have $x_0=(-y,0)$ and
  $x_+=(0,c)$ where $0<y\leq C\delta t_+$. We verify that
  \begin{displaymath}
    \tilde\vp_{x_0}(t) - x_0 - t v_0 = \int_0^t\int_0^s 
    \partial^2_t{{\tilde{\varphi}_x}}(u)duds.
  \end{displaymath}
    $\tilde\varphi_{x_0}(t)$ will be outside
  of the convex domain $X$ as soon as its first component becomes
  non-negative, which implies that
  \begin{displaymath}
       -y + t v\cdot e_1 + e_1\cdot \int_0^t\int_0^s 
    \partial^2_t{{\tilde{\varphi}_x}}(u)duds \leq 0,
  \end{displaymath}
  or equivalently
  \begin{displaymath}
    t_+ t \leq C_3 t_+ \delta + C_4 t^2. 
  \end{displaymath}
  Now take $t=2C_3\delta$ so that the above constraint becomes
  \begin{displaymath}
    t_+ \leq 4C_4 C_3 \delta. 
  \end{displaymath}
  It remains to choose $t_+ > 4C_4 C_3 \delta$ to obtain the existence
  of a time $t$ so that $\tilde\varphi_{x_0}(2C_3\delta)$ is outside
  of $X$. This shows that the distance traveled by $\tilde\varphi_{x_0}$
  is of size $\delta$ so that $|x_+(x)-\tilde x_+(x)|\leq C\delta$.
  
  We have treated the case $t_+(x)>\alpha \delta$ for some constant
  $\alpha$ sufficiently large. It remains to address the case
  $t_+(x)<\alpha \delta$. For this, a sufficiently positive curvature
  of $\pa X$ is necessary and we use Lemma \ref{lem:propxt}. Indeed,
  we know that
  \begin{displaymath}
    |\varphi_x(t_+(x))-\tilde\varphi_x(t_+(x))| 
    \leq C\alpha \delta^2.
  \end{displaymath}
  We also know that $\beta(x_+(x))\cdot n(x_+(x)) \geq 2\alpha C_1
  \delta$. We can then invoke the second result of Lemma
  \ref{lem:propxt} and obtain that $|x_+(x)-\tilde x_+(x)|\leq
  C\delta$.
  
  At this stage, we have thus proved that independent of $t_+(x)$,
  $|x_+(x)-\tilde x_+(x)|\leq C\delta$ with $C$ a constant independent
  of $\delta$. The proof of the above result shows that
  $|t_+(x)-\tilde t_+(x)|\leq C\delta$ as well.  Higher order
  derivatives are now treated in a similar fashion. We have seen that 
  \begin{displaymath}
    |W(t_+)-\tilde W(t_+)|\leq C t_+\|D_x\beta-D_x\tilde\beta\|_{C^0(\bar X)}.
  \end{displaymath}
  Since $\tilde W(t)$ is of class $C^1$ and $|t_+(x)-\tilde
  t_+(x)|\leq C\delta$, we deduce that 
  \begin{displaymath}
    |W(t_+)-\tilde W(\tilde t_+)| \leq C \|\beta-\tilde\beta\|_{C^1(\bar X)}
  \end{displaymath}
  which is equivalent to
  \begin{displaymath}
    \|D_x x_+ - D_x \tilde x_+\|_{C^0(\bar X)} 
   \leq C \|\beta-\tilde\beta\|_{C^1(\bar X)}.
  \end{displaymath}
  Higher-order derivatives are treated in exactly the same manner providing
  a bound for $k-1$ derivatives of $x_+-\tilde x_+$ in the uniform norm. 

  The error on $t_+$ is obtained as follows. We note that
  \begin{displaymath}
    \vp_x(t_+)-\tilde\vp_x(\tilde t_+) = (x_+-\tilde x_+)(x).
  \end{displaymath}
  After differentiation in space, we obtain
  \begin{displaymath}
    W(t_+(x)) D_x t_+(x) - \tilde W(\tilde t_+(x)) D_x \tilde t_+(x)
    = D_x (x_+-\tilde x_+)(x).
  \end{displaymath}
  Since $W$ is Lipschitz, and $t_+-\tilde t_+$ is small, this implies
  that
  \begin{displaymath}
    W(t_+(x)) (D_xt_+(x)-D_x\tilde t_+(x)) = O(\delta).
  \end{displaymath}
  We have that $\dot W$ is of order $|\k|^{-1}$ and that $W(0)=I$ so
  that for $|\k|$ sufficiently large, $W(t_+(x))$ is invertible. This
  implies
  \begin{displaymath}
    (D_xt_+(x)-D_x\tilde t_+(x)) = O(\delta).
  \end{displaymath}
  Higher-order derivative are treated in the same manner by using the
  Leibniz product rule and the invertibility of $W(t_+(x))$. This
  concludes the proof of the lemma.
\end{proofof}

\begin{proofof} {\it Lemma \ref{lem:propxt}}.
  Instead of running the characteristics forward from $x$ to $x_+(x)$,
  we run the characteristics backwards from $x_+(x)$ to $x_-(x)$ and
  show that the time spent in $X$ is controlled by the angle the
  trajectories makes with the normal to $X$ at $x_+(x)\in\pa X$. More
  precisely, we set $y=x_+(x)$ and $v=\dvp_x(t_+(x))=\beta(x_+(x))$
  and run characteristics backwards.

From the equality
\begin{displaymath}
  \dvp_x(-t) = v - \int_0^t\ddot\vp_x(-s)ds,
\end{displaymath}
we deduce that 
\begin{displaymath}
  |\dvp_x(-t)-v|\leq \dfrac{C_0t}{|\k|},
\end{displaymath}
and hence
\begin{displaymath}
  \big|\vp_x(-t)-\big(y-tv\big)\big|\leq \dfrac{C_0t^2}{2|\k|}.
\end{displaymath}
Let $t_m$ the time it takes
from $x_+(x)$ to $x_-(x)$. We obviously have that $t_+(x)\leq t_m$.
Let $B_y(R)$ the (unique) ball of radius $R$ tangent to $\pa X$ at
$y\in\pa X$ and such that $X\subset B_y(R)$.

In a system of coordinates with $B_y(R)$ centered at $0$ and $v=e_1$,
we find that $|y-tv|^2=R^2+|v|^2t^2-2rt v\cdot n(y)$. We want
$\vp_x(-t)\in X\subset B_y(R)$. This imposes that
\begin{displaymath}
  |y-tv| \leq R + \dfrac{C_0}{2|\k|} t^2.
\end{displaymath}
Let us define $t_M$ as the maximal time a trajectory spends in $X$,
which is a bounded quantity. Then the above imposes that
\begin{displaymath}
  \begin{array}{l}
     |v|^2 t^2-2Rt v\cdot n(y) \leq \dfrac{RC_0}{|\k|} t^2
     + \dfrac {C_0^2 t_M^2}{4|\k|^2} t^2 \leq |v|^2\rho t^2,
  \end{array}
\end{displaymath}
and in other words that $t_m$ is bounded by
$\frac{2R}{|v|^2(1-\rho)}v\cdot n(y)$.

For the second result, we define $B_x(R)$ as the ball of radius $R$
tangent to $\pa X$ at $x$ and such that $X\subset B_x(R)$. We again
have that $|\varphi_{x_0}(t)-(x_0+v_0t)|$ is bounded by $\rho|v_0|^2
t^2$.  In a system of coordinates where $n(x)=e_2$, $x=Re_2$, and
$v_0=(v_0\cdot e_1)e_1+(v_0\cdot e_2)e_2$, we obtain that
$\varphi_{x_0}(t)\in X\subset B_x(R)$ implies that
\begin{displaymath}
  |x+tv_0|\leq |x-x_0| + R + \dfrac{C_0}{2|\k|}t^2.
\end{displaymath}
This is equivalent to
\begin{displaymath}
  2tRv_0\cdot e_2 + t^2 |v_0|^2 \leq |x-x_0|^2+\dfrac{C_0^2t_M^2}{4|\k|^2}t^2
   +2|x-x_0|R+2|x-x_0|\dfrac{C_0t_M^2}{2|\k|} +\dfrac{RC_0}{|\k|}t^2.
\end{displaymath}
This implies that
\begin{displaymath}
  (1-\rho)|v|^2t^2 \leq C|x-x_0| + 2t R |v\cdot e_2|
   \leq C(\delta^2+t\delta),
\end{displaymath}
for some constants $C$ that can be made explicit. Solving this
quadratic inequality yields that $t_m$ is bounded as prescribed.
\end{proofof}

%
\subsection{ISID with two real-valued measurements}
\label{sec:isid2}
%

We are now in a position to state our main stability result in the
presence of one well chosen complex-valued internal data.  We fix $\k$
and let $\rho=\k+i\k^\perp$. We define $u$ as the solution to
\eqref{eq:u} with the complex-valued illumination $g$ close to
$u_{\rho|\pa X}$.  The complex-valued internal data are then $d=\mu
u$.  As before, we assume that $k\geq1$ and $\mu\in C^{k+1}(\bar X)$
so that $d\in C^{k+1}(\bar X;\Cm)$, which implies that $\beta$ and
$\gamma$ defined in \eqref{eq:betagamma0} are of class $C^k(\bar X)$.

The results of Proposition \ref{prop:regmu} yield the following result:
\begin{theorem}
  \label{thm:reg2m}
  Let us assume that $(\mu,q)$ and $(\tilde\mu,\tilde q)$ are elements
  in $\mathcal P$ and that $\|g-u_{\rho|\pa X}\|_{C^0(\bar X)}$ is
  sufficiently small so that $u$ does not vanish on $X$.  Under the
  hypotheses of Proposition \ref{prop:regmu} and assuming that
  $|\k|\geq|\k_0|$ with $|\k_0|$ sufficiently large, we have that
  \begin{equation}
    \label{eq:regmu2m}
    \|\mu-\tilde\mu\|_{C^{k-1}(\bar X)}\leq 
    C \|d-\tilde d\|_{(C^k(\bar X;\Cm))}.
  \end{equation} 
  Moreover, we have the following stability result provided that
  $k\geq3$:
  \begin{equation}
    \label{eq:stabqmu}
    \|q-\tilde q\|_{C^{k-3}(\bar X)} 
    \leq  C \|d-\tilde d\|_{(C^k(\bar X;\Cm))}.
  \end{equation}
\end{theorem}
\begin{proof}
  Let us define $\mu_0=d_{|\pa X}/g$ and $\tilde\mu=\tilde d_{|\pa
    X}/g$. By assumption, $g$ does not vanish on $\pa X$. We thus
  deduce that $\|\mu_0-\tilde\mu_0\|_{C^k(\pa X)}$ is controlled by
  $\|d-\tilde d\|_{(C^k(\bar X;\Cm))}$. The rest of the proof of the
  theorem is a direct consequence of the results obtained in the
  preceding section and of the proof of Theorem \ref{thm:stab2n}.
\end{proof}

We now make a few comments on the differences between the $2n-$ and
$2-$ internal data settings. In the presence of $n$ fields, the
geometry of $X$ is allowed to be rather general. In contrast, the $2-$
data setting requires much stronger convexity assumptions on $X$ to
avoid that integral curves of the vector field be too close to the
boundary $\pa X$ for too long, which would result in a severe lack of
stability.  Since the integral curves of the vector fields $\beta$ and
$\tilde\beta$ are not known a priori, more (Lipschitz) regularity is
required on the vector fields to ensure that information propagates
along near-by trajectories. This is the reason for the replacement of
$k$ in Theorem \ref{thm:stab2n} by ``$k-1$'' in Theorem
\ref{thm:reg2m}.

\section{Inverse Diffusion with Internal Data (IDID)}
\label{sec:idid}

We now return to the diffusion equation with unknown diffusion
coefficient $D$ and unknown absorption coefficient $\sigma_a$:
\begin{displaymath}
  -\nabla \cdot D \nabla u + \sigma_a u =0, \quad X, \qquad
   u = g,\quad \pa X.
\end{displaymath}
The theorems stated in section \ref{sec:main} are straightforward
consequences of the results presented in this section in a slightly
more general setting.

Using the standard Liouville change of variables, $v=\sqrt D u$ solves
\begin{displaymath}
  \Delta v + qv =0,
\end{displaymath}
with
\begin{displaymath}
  q=-\dfrac{\Delta\sqrt D}{\sqrt D} -\dfrac{\sigma_a}{D}.
\end{displaymath}
The internal data in photoacoustics are given by
\begin{displaymath}
  d=\sigma_a u = \dfrac{\sigma_a}{\sqrt D} v = \mu v,\qquad
  \mu :=  \dfrac{\sigma_a}{\sqrt D}.
\end{displaymath}
We assume we know $\sqrt D$ on $\partial X$. This allows us to prescribe
$v$ on $\partial X$ and thus to reconstruct $\mu$ and $q$ as in the 
preceding section. Then we find that 
\begin{equation}\label{eq:recDfromq}
  -\Delta \sqrt D - q \sqrt D =\mu,
\end{equation}
so we can solve for $\sqrt D$ and then get $\sigma_a=\mu\sqrt D$.


Let us recall that $\sqrt D\in Y=H^{\frac n2+k+2+\eps}(X) \subset
C^{k+1}(\bar X)$.  We assume that $(D,\sigma_a)\in {\mathcal M}$ for
some $M$.  This implies that $(q,\mu)\in {\mathcal P}$ for some $P$.
Indeed, that $\sqrt D$ solves \eqref{eq:recDfromq} implies that
$\Delta+q$ whose inverse is compact does not have $0$ as an
eigenvalue. Here, as always, $k\geq1$.

The above calculations show the unique reconstruction of
$(D,\sigma_a)$ from internal data for well-chosen boundary
distributions as stated in Theorem \ref{thm:uniq000}.  From Theorem
\ref{thm:stab2n} in the $2n-$internal data setting, we get the
following result.
\begin{theorem}
  \label{thm:stabDsigma}
  Let $k\geq2$ and assume that $(D,\sigma_a)$ and $({\tilde D},\tilde
  \sigma_a)$ are in $\mathcal M$ with $D_{|\pa X}=\tilde D_{|\pa X}$
  on $\pa X$.  Then there is an open set of $2n$ real valued boundary
  values $g$ in $C^{k,\alpha}(\pa X)$ for $\alpha>\frac12$ such that
  we have the stability estimate
  \begin{equation}
  \label{eq:stabDsigma2}
   \|D-\tilde D\|_{C^{k}(X)} + \|\sigma_a-\tilde\sigma_a\|_{C^{k}(X)}
   \leq C(M,k) \|d-\tilde d\|_{(C^{k+1}(X))^{2n}}.
  \end{equation}
\end{theorem}
\begin{proof}
  The main result consists of getting the stability on $D$ mentioned
  above. Since $k\geq2$, we have stability of the reconstruction of
  $q$ in $C^{k-2}(\bar X)$ and of $\mu$ in $C^k(\bar X)$ provided that
  the boundary conditions are well-chosen. Thus we have
  \begin{displaymath}
    -(\Delta +q) (\sqrt{D}-\sqrt{\tilde D}) = \mu-\tilde \mu
   + (q-\tilde q)\sqrt{\tilde D}.
  \end{displaymath}
  By elliptic regularity, we deduce that $(\sqrt{D}-\sqrt{\tilde D})$
  is bounded in $C^{k}(\bar X)$, and hence the result. 
\end{proof}


Finally, we deduce from Theorem \ref{thm:reg2m} the following result
\begin{theorem}
  \label{thm:reg2mD}
  Under the hypotheses of Theorem \ref{thm:reg2m} and those of Theorem
  \ref{thm:stabDsigma}, we obtain in the 2-internal data setting the
  following result.  Let $k\geq3$ and $(D,\sigma_a)$ and $(\tilde
  D,\tilde\sigma_a)$ be in $\mathcal M$. Let us assume that $D_{|\pa
    X}=\tilde D_{|\pa X}$.
  
  Then there is an open set of $2$ real-valued boundary conditions $g$
  in $C^{k,\alpha}(\pa X)$ for $\alpha>\frac12$ such that we have the
  stability estimate
  \begin{equation}
  \label{eq:stabDsigma}
   \|D-\tilde D\|_{C^{k-1}(X)} + \|\sigma_a-\tilde\sigma_a\|_{C^{k-1}(X)}
   \leq C(M,k) \|d-\tilde d\|_{C^{k}(X;\Cm)}.
  \end{equation}
\end{theorem}
The proof of the above theorem is the same as that of Theorem
\ref{thm:stabDsigma}.

\section*{Acknowledgment} 
The authors would like to thank John Sylvester for interesting
discussions on the inverse diffusion problem with internal data.  GB
was supported in part by NSF Grants DMS-0554097 and DMS-0804696.  GU
was supported in part by NSF and a walker Family Endowed
Professorship.


\end{document}